\documentclass[12pt]{article}
\usepackage[centertags]{amsmath}
\usepackage{amsfonts}
\usepackage{amssymb}
\usepackage{amsthm}
\usepackage{amsmath,mathrsfs}
\usepackage{amsmath,amscd}
\usepackage{mathtools}
\usepackage[matrix,arrow,curve,cmtip]{xy}
\title{\textbf{Homotopy Wave Function in Algebraic Geometry}}
\author{Renaud Gauthier \footnote{2020 Math. Subj. Class: 18N60, 18F20, 55P55, 33E15, 53E40 . Keywords: wave function, Segal topos, derived stack, homotopy shape} \\ \\}
\theoremstyle{definition}

\newcommand{\beq}{\begin{equation}}
\newcommand{\eeq}{\end{equation}}

\newcommand{\rarr}{\rightarrow}

\newcommand{\xrarr}{\xrightarrow}

\newcommand{\cA}{\mathcal{A}}

\newcommand{\cC}{\mathcal{C}}

\newcommand{\cX}{\mathcal{X}}

\newcommand{\bR}{\mathbb{R}}

\newcommand{\Hom}{\text{Hom}}
\newcommand{\Ho}{\text{Ho}\,}

\newcommand{\Loc}{\text{Loc}}

\newcommand{\Map}{\text{Map}}

\newcommand{\op}{\text{op}}

\newcommand{\Set}{\text{Set}}

\newcommand{\Top}{\text{Top}}
\newcommand{\uHom}{\underline{\Hom}}

\newcommand{\AffC}{\cA \text{ff}_{\cC}}

\newcommand{\Comm}{\text{Comm}}

\newcommand{\CommC}{\text{Comm}(\cC)}

\newcommand{\dAff}{\text{dAff}}
\newcommand{\dAffC}{\dAff_{\cC}}
\newcommand{\dAffCtildetau}{\dAff_{\cC}^{\;\sim \, , \, \tau}}

\newcommand{\dkAff}{\text{d}k\text{-Aff}}
\newcommand{\dkAfftildetau}{\dkAff^{\,\sim \, , \, \tau}}
\newcommand{\dSt}{\text{dSt}}

\newcommand{\RHom}{\mathbb{R} \uHom}

\newcommand{\RHomgeom}{\RHom^{\text{geom}}}

\newcommand{\RuHom}{\bR \uHom}

\newcommand{\sPr}{\text{sPr}}

\newcommand{\SetD}{\Set_{\Delta}}

\newcommand{\skMod}{\text{s}k\text{-Mod}}
\newcommand{\skCAlg}{\text{s}k\text{-CAlg}}

\newcommand{\SeT}{\text{SeT}}

\newcommand{\SeCat}{\text{SeCat}}

\newcommand{\XF}{\cX_{F/}}
\newcommand{\HXF}{H_{\XF}}

\begin{document}
\maketitle
\begin{abstract}
We propose the homotopy shape of the Segal topos of derived stacks over simplicial $k$-algebras as the higher homotopical generalization of the concept of wave function in Quantum Mechanics.
\end{abstract}

\section{Introduction}
It is fairly well-known that there is such a thing as a wave-particle duality in Quantum Mechanics, according to which particles can equivalently be described as being corpuscular, or be represented by waves functions. Recent experiments, \cite{LFRFG} being just one example, show that individual particles exhibit wave-like behaviors, which renewed this interest in the interpretation of Quantum Mechanics, and about the wave-particle duality in particular. Among the many questions that surround the wave-particle interpretation of Quantum Mechanics, one problem that seems not to be settled is that of the ``wave function collapse", according to which when we do not measure the position of a particle, it is in its wavy state, provided by a wave function, and once we observe that particle, it turns into just that, a particle, and the wave describing that particle is gone, and its wave function along with it. In the present paper we provide an explicit construction of a higher version of the concept of wave function that addresses those points.\\

The wave function $\psi$ in Quantum Mechanics provides all the information there is to know about a given particle of interest, and its norm squared provides the probability density of observing that particle in one state or another. In other terms, the wave function fully describes a phenomenon, and this in a smooth fashion since its norm squared is an integrable distribution. On our part we have developed a homotopical model of natural phenomena using derived stacks over simplicial $k$-algebras for a commutative ring $k$, the collection of which forms a Segal topos $\cX$. Given one phenomenon $F \in \cX$, we argue the comma Segal topos $\XF$ provides the categorical manifestation of $F$ within $\cX$. In other terms it presents $F$ within $\cX$ as a dynamical phenomenon yet to become something else, the objects $F \rarr G$ of the comma category providing those future states $G$ the stack $F$ will turn into. Thus $\XF$ fully describes $F$ dynamically within its ambient universe $\cX$, in a pedestrian sense. This being done, next we want to exhibit the homotopical character of such a representation of $F$. We argue one object that provides such an enhanced view of $F$ is the homotopy shape $\HXF$ of $\XF$. Recall how this is constructed (\cite{SeT}): for $X$ a Segal topos, $Y$ a CW-complex:
\beq
H_{X}(Y) := |\RHomgeom(X,\Loc(Y))| \nonumber
\eeq
is the geometric realization (read diagonal simplicial set) of the Segal category of geometric morphism between $X$ and the Segal topos $\Loc(Y)$ of locally constant stacks on $Y$. Why would this provide an acceptable homotopical presentation of $X = \HXF$ in our case? For one thing it would have to be a homotopical object, which it is, but most importantly by being ``acceptable" we mean to say it would have to capture the slightest variations in $F$ encoded in $\XF$ for the sake of completeness, and those are best measured relative to objects that are stable, such as locally constant stacks, thus if there is an object that captures any smooth dynamics it is the homotopy shape. Collecting things together, $\XF$ presents $F$ dynamically in full in $\cX$, and its homotopy shape $\HXF$ displays this homotopically. The functor $\HXF: \Top \rarr \Top$, $\Top$ being the Segal category of simplicial sets, is represented by an associated object $\HXF \in \RuHom(\Top,\Top)$ that we refer to as the \textbf{homotopy wave function} for the phenomenon $F$. In simple cases, this at least contains the information given by a wave function in Quantum mechanics, which justifies the appellation we gave to $\HXF$.\\

In a first time we remind the reader how is $\cX = \dSt(k)$ constructed and why would that Segal topos correspond to a modelization of natural phenomena. At the same time we will show that this is a homotopical object as well. Given one phenomenon $F \in \cX$, we will focus on the Segal topos $\XF$, which therefore describes $F$ in full dynamically, but also homotopically. However, the homotopical character of $F$ within $\cX$ is best identified once it is extracted using the homotopy shape $\HXF$, and this we will discuss last. Finally we argue the wave function in Quantum Mechanics should at least be encoded in this homotopy shape, which therefore provides us with a higher homotopical wave function.\\

\textbf{Relation to other work}: in \cite{K} a quantum ensemble is defined, which if ``coherent" enough not only defines an implied notion of smoothness but can also be described by a wave function as it is generally understood in Quantum Mechanics. No collapse is invoked, since the quantum statistical ensemble always exist. What we present is a higher homotopical construction of not only the background setting for the existence of a wave function, but of the actual homotopical wave function itself, which in simple cases boils down to a formalism similar to, but not identical with, the one alluded to in \cite{K}.

\section{Segal topos of derived stacks $\cX = \dSt(k)$}
\subsection{Construction}
In \cite{RG} we argued that natural laws can be modeled by simplicial algebras over a commutative ring $k$. We chose to work with simplicial objects since it is natural to embed classical notions in a derived setting to avoid having to deal with singularities. To manifest those laws in nature we need coherent manifestations, and for that purpose stacks are all indicated. Thus starting from the model category $\skCAlg = \Comm(\skMod)$ of simplicial algebras over $k$, we first consider prestacks $F: \skCAlg \rarr \SetD$. To account for the fact that in nature some algebraic equivalences are not necessarily observable, we take a simplicial localization of such model categories, which throws us into the realm of Segal categories (\cite{SeT}); we work with $L (\skCAlg)^{\op} = \dkAff$, the Segal category of derived affine stacks (\cite{HidSt}, \cite{AGmodCat}, \cite{HAGI}, \cite{HAGII}), and $\Top = L \SetD$ the Segal category of simplicial sets, and we therefore consider prestacks $F \in \RuHom(\dkAff^{\op}, \Top) = \widehat{\dkAff}$. If we put a Segal topology on $\dkAff$, we obtain a Segal site, hence a Segal category of stacks $\dkAfftildetau$, which turns out to be a Segal topos usually denoted $\dSt(k):=\cX$.

\subsection{Homotopical behavior of $F \in \cX$}
Models for derived stacks $F \in \cX$ are given by functors $F: \skCAlg \rarr \SetD$ that preserve equivalences, so if there is but a tiny smooth variation in a natural law $A \in \skCAlg$, then $FA$ changes smoothly as well. Now observe that this is just a response of a given phenomenon $F$ to a variation in $A$. The complete picture $F$ provides is really given by its essential image, so one should consider all $FA$ for $A$ running over $\skCAlg$. As a matter of fact this is the right object to consider, the reason being that if $A$ is altered in any way, $\skCAlg$ being a category, other algebras within will be modified as well, so one should consider all their associated images $FA$ simultaneously. There is actually a way to do just this: take the join (\cite{J}) of all such simplicial sets, so that we consider:
\beq
\sum_{A \in \skCAlg} FA := *_{A \in \skCAlg} FA \in \SetD \nonumber
\eeq
Since $F$ preserves equivalences, so does $* F$ by definition. Now in what sense would such equivalences correspond to smooth deformations in the pedestrian sense? Recall that equivalences in $\SetD$ are defined to be weak equivalences in the category of topological spaces, in which category maps are continuous, so if $f:X \rarr Y$ is a continuous map of topological spaces, we say this is a weak equivalence if for all $n\geq 0$, for any base point $x \in X$, we have isomorphisms $\pi_n(X,x) \xrarr{f_*} \pi_n(Y,fx)$, and each higher homotopy group $\pi_n(X,x)$ is the set of pointed homotopy classes of pointed continuous maps $(S^n,*) \rarr (X,x)$, and homotopies are continuous maps. What this really means is that equivalences in $\SetD$ are smooth, read are defined from homotopies. Thus to summarize if algebras are deformed smoothly, that directly translates into $\sum F = * F$, the materialization of the phenomenon $F$, to also change smoothly.

\subsection{Smooth flows within $\cX$}
We now discuss flows in $\cX$, which we recall are generated by composing objects $\psi \in \RuHom(\cX,\cX) = \delta \cX$ (\cite{RG4}), q-deformations to be precise. We can regard those as tangent vectors that generate smooth deformations of objects $F \in \cX$, and talking about a flow in $\cX$ through $F$ this means considering all those successive actions on $\cX$ by objects of $\delta \cX$, resulting in an integral curve that we refer to as the flow of $F$. To argue in favor of having a smooth flow, consider objects $\psi \in \delta \cX$ acting on $F \in \cX$ as $F \mapsto \psi F$. If those objects $\psi$ are small enough (q-deformations), such maps $F \mapsto \psi F$ are infinitesimal, and in particular can be considered to also include equivalences in $\cX \in \SeCat$. Since it is the flows we consider to be smooth, it is really $\psi \in \delta \cX = \RuHom(\cX,\cX)$ as a map of Segal categories we should focus on. Assuming $\psi$ to be an equivalence, recall how those are defined (\cite{SeT}): a morphism $f:A \rarr B$ between Segal categories is an equivalence if it is fully faithful, that is for all objects $a,a'$ of $A_0$, we have an isomorphism $A_{(a,a')} \rarr B_{(fa,fa')}$ in $\Ho(\SetD)$, or equivalently an equivalence in $\SetD$, and $f$ is also essentially surjective, that is it is essentially surjective at the level of homotopy categories. Thus $\psi$ being an equivalence as a map can be regarded as being smooth, and composites $\Phi$ of such objects are again smooth.\\

Now this produces a deformation in $\cX$, which recall parametrizes phenomena, or equivalently said is a moduli space of natural phenomena. However, those deformations do not occur in $\SetD$, where natural phenomena actually take place, they really happen on $\cX$. To make the transition to the natural realm, observe that flows $\Phi F$ through $F$ for $\Phi \in \delta \cX$ are defined objectwise, so the corresponding materialization in $\SetD$ of such a flow is given by $\Phi FA$ for all $A \in \skCAlg$, and taking the join of all such simplicial sets we see that $* \Phi F$ provides a smooth deformation of $F$ as a whole.\\

Regarding those flows, any object $\Phi \in \delta \cX$ will generate a corresponding deformation of $\cX$, and consequently will alter a given phenomenon $F$ within $\cX$ in a certain fashion, whether this be just incremental or a smooth flow. However, given the past dynamics within $\cX$, some future flows will be more likely than others, and this is important to keep in mind. In other terms the present state of affairs is given by $\cX$ as a whole, future dynamics therein is encoded in $\delta \cX$, but consequently carries a probabilistic character, since some objects $\Phi$ within that generate flows are more likely than others. On another note, one does not have to consider $\delta \cX$ as some additional data, since we can naturally identify $\cX$ with $id_{\cX} \in \RuHom(\cX,\cX) = \delta \cX$, thus seeing $\cX$ as an object of $\delta \cX$ those possible flows on $\cX$ are part of the ambient picture $\cX$ finds itself in, namely $\delta \cX$.

\section{Upgrade to $\XF$}
Focusing on one particular phenomenon $F$ within the ambient Segal topos $\cX$, we argue future dynamics thereof is provided by objects of $\XF$ (since those are of the form $F \rarr G$), which still encapsulates all there is to know about $F$, as well as its homotopical character as argued above, whether it be relative to smooth deformations of natural laws in $\skCAlg$, or that which results from flows in $\cX$ through $F$ generated by objects $\Phi \in \delta \cX$. These two types of deformations provide us with the future dynamics of $F$, that comprises moves that are more or less likely as pointed out above, so saying $\XF$ is a homotopical object directly translates into working with a geometric object that is inherently probabilistic in nature. Since $\XF$ describes $F$ in full, and describes its probabilistic behavior, we tentatively regard it  as a homotopical wave function for $F$ in the rough.

\section{Homotopical extraction with the Homotopy Shape $\HXF$}
As a further upgrade, we argue we have to make the homotopical character of $\XF$ manifest and there is no better way to achieve this than by comparing the Segal topos $\XF$ to other Segal topos that are as inert as one would like to imagine, since from the perspective of such topos, the slightest deformation can easily be picked up. From \cite{Gal}, \cite{Hom_HiCat} and \cite{SeT}, recall that if $X$ is a CW-complex, $\sPr(X)$ is the simplicial category of simplicial presheaves on $X$ with its local projective model structure, one can consider $\text{PrLoc}(X) \subset \sPr(X)$ the full subcategory of locally constant presheaves on $X$ and then consider its associated Segal category $\Loc(X) = L(\text{PrLoc}(X))$. As shown in \cite{SeT}, not only is that a Segal topos, but we have a fully faithful embedding $\Loc: \Ho(\Top) \rarr \Ho(\Top)$, $X \mapsto \Loc(X)$. Further in the same reference it is shown that for two CW-complexes $X$ and $Y$, $\Map(X,Y) \cong |\RHomgeom(\Loc(X),\Loc(Y))|$ in $\Ho(\SetD)$, thus the correct object of comparison between two CW-complexes at the level of Segal topos, an ambient environment within which they embed faithfully, is the geometric realization $|\RHomgeom(\Loc(X),\Loc(Y))|$, which is none other than a geometric shape. Thus we first consider $\RHomgeom(\XF,\Loc(Y))$ for all possible CW-complexes $Y$ for the sake of comparing $\XF$ to locally constant objects in $\SeT$, the Segal category of Segal topos, and we actually take the geometric realization of such a Segal groupoid, giving rise to a functor $\HXF: \Top \rarr \Top$ which to any CW-complex $Y$ associates $\HXF(Y) = |\RHomgeom(\XF,\Loc(Y))|$, which is referred to as the homotopy shape of $\XF$. Finally we take the corresponding object of $\RuHom(\Top,\Top)$ that we still denote with the same notation. This is our higher homotopical wave function. To be fair, the homotopy shape $H_X$ for $X \in \SeT$ is well-defined if in addition $X$ is t-complete, which is the case for the Segal topos of stacks such as the one we are using.\\

The argument that this provides the classical wave function of Quantum Mechanics is fairly simple; the homotopy shape encodes all there is to know about the phenomenon $F$, homotopical behavior included, and we argued this is equivalent to regarding that as probabilistic information. Since $\psi_F$ in Quantum Mechanics plays the same role, but without proof of the contrary does not contain higher homotopical information, it must at least be contained within $\HXF$ information-wise. In other terms our homotopy shape reduces to the classical wave function in Physics, or equivalently, the generalization to higher mathematics of the wave function as we know it in Quantum Mechanics is provided by the homotopy shape of the phenomenon the said wave function is supposed to describe.\\

What of the ``wave-function collapse"? We argue this is just measuring a slice in the Segal topos $\XF$, and the latter always exist, regardless of any measurement. In other terms if there is collapse it is only apparent, we measure what we expect to observe, but we fail to observe a flow simply because at any moment we observe only single manifestations.

\section{Generalizations}
We focused exclusively on $\dSt(k)$, with base model category $\skMod$, with $\skCAlg = \Comm(\skMod)$, the idea being that simplicial algebras model natural laws. One could rightly argue that this is far from sufficient to model natural phenomena in their integrality. This can be remedied by using simplicial monoidal model categories $\cC$ satisfying the conditions of a homotopical algebraic concept according to \cite{HAGII} for the sake of having derived stacks defined on them, but one thing we still keep from the above picture is that the materialization of such foundational behavior given by $\cC$ is provided by derived stacks. Thus if $\cC$ is such a base model category, we consider its category $\Comm(\cC)$ of commutative monoids, which is again a model category since $\cC$ satisfies the conditions of a HA context. We consider its opposite category $\CommC^{\op} = \AffC$, take a simplicial localization thereof (or a Segal localization) $L\AffC = \dAffC$, on which we put a Segal topology, that is a Grothendieck topology on its homotopy category, which gives rise to a Segal category $\dAffCtildetau$ of derived stacks, modeled by functors $F: \Comm(\cC) \rarr \SetD$ and we proceed from there. All the work above carries over in this general case.\\

Now suppose a given natural phenomenon $\varphi$ depends on natural laws as well as other information other than natural laws. Regard $\varphi$ as a functorial manifestation of $\skCAlg$, as well as other simplicial symmetric monoidal model categories $\cC_1,\cdots, \cC_N$, and each of those categories $\cC_i$, $1 \leq i \leq N$ gives rise to a Segal category of stacks $\text{dAff}_{\cC_i}^{\, \sim \, , \, \tau_i}:= \cX_i \in \SeT$ for some appropriate Segal topology $\tau_i$ on $\dAff_{\cC_i}$, and within $\cX_i$, there is a derived stack $F_i$ describing the $\cC_i$-aspect of $\varphi$, or equivalently said, $\varphi$ is the manifestation of $F_i \in \cX_i$ in the natural realm. For each of those Segal topos we consider the comma category $\cX_{i,F_i /} \in \SeT$, we consider their product $\cX_{\varphi/}: = \prod_i \cX_{i,F_i/} \in \SeT$, which effectively presents the future dynamics of $\varphi$ from the perspective of its different constituting stacks $F_i$ within their respective topos $\cX_i$. Finally we consider the homotopy shape $H_{\cX_{\varphi/}} \in \RuHom(\Top,\Top)$, which we regard as the homotopy wave function of $\varphi$, a phenomenon that depended on natural laws and possibly other base model categories.

\bigskip
\footnotesize
\noindent
\textit{e-mail address}: \texttt{rg.mathematics@gmail.com}.

\end{document}